\documentclass[12pt,A4paper]{amsart}
\usepackage{amsmath,amsthm,amssymb,xcolor}
\usepackage[bookmarksnumbered,colorlinks]{hyperref}
\hypersetup{colorlinks=true,linkcolor=red,anchorcolor=green,citecolor=cyan}

\newtheorem{theorem}{Theorem}[section]
\newtheorem{lemma}[theorem]{Lemma}
\theoremstyle{remark}

\numberwithin{equation}{section}
\newcommand{\N}{\mathbb{N}}

\newcommand{\C}{\mathbb{C}}
\newcommand{\D}{\mathbb{D}}
\newcommand{\T}{\mathbb{T}}

\begin{document}

\title[\tiny Spaceability of special families of null sequences of holomorphic functions]{Spaceability of special families of null sequences of holomorphic functions}

\author[Bernal]{L.~Bernal-Gonz\'alez}
\address[L. Bernal-Gonz\'alez]{\mbox{}\newline \indent Departamento de An\'alisis Matem\'atico \newline \indent Facultad de Matem\'aticas
\newline \indent Instituto de Matem\'aticas de la Universidad de Sevilla (IMUS)
\newline \indent Universidad de Sevilla
\newline \indent Avenida Reina Mercedes s/n, 41012-Sevilla (Spain)}
\email{lbernal@us.es}

\author[Calder\'on]{M.C.~Calder\'on-Moreno}
\address[M.C.~Calder\'on-Moreno]{\mbox{}\newline \indent Departamento de An\'alisis Matem\'atico \newline \indent Facultad de Matem\'aticas
\newline \indent Instituto de Matem\'aticas de la Universidad de Sevilla (IMUS)
\newline \indent Universidad de Sevilla
\newline \indent Avenida Reina Mercedes s/n, 41012-Sevilla (Spain)}
\email{mccm@us.es}

\author[L\'opez-Salazar]{J.~L\'opez-Salazar}
\address[J.~L\'opez-Salazar]{\mbox{}\newline \indent Departamento de Matem\'atica Aplicada a las Tecnolog\'ias  \newline \indent de la Informaci\'on y de las Comunicaciones
\newline \indent Escuela T\'ecnica Superior de Ingenier\'ia y Sistemas de Telecomunicaci\'on
\newline \indent Universidad Polit\'ecnica de Madrid
\newline \indent Nikola Tesla s/n, 28031-Madrid (Spain)}
\email{jeronimo.lopezsalazar@upm.es}

\author[Prado]{J.A.~Prado-Bassas}
\address[J.A.~Prado-Bassas]{\mbox{}\newline \indent Departamento de An\'alisis Matem\'atico
\newline \indent Facultad de Matem\'aticas
\newline \indent Instituto de Matem\'aticas de la Universidad de Sevilla (IMUS)
\newline \indent Universidad de Sevilla
\newline \indent Avenida Reina Mercedes s/n, 41012-Sevilla (Spain)}
\email{bassas@us.es}



\begin{abstract}
In this note, we consider the space \,$H(\Omega)^{\N}$ \,of sequences of holomorphic functions on an open set $\Omega\subset \C$. If \,$H(\Omega)$ \,is endowed with its natural topology and \,$H(\Omega)^{\N}$ \,is endowed with the product topology, then it is proved the existence of
two closed infinite dimensional vector subspaces of \,$H(\Omega)^{\N}$ \,such that all nonzero members of the first subspace are sequences tending to zero pointwisely but not compactly on \,$\Omega$ and all nonzero members of the second subspace are sequences tending to zero
compactly but not uniformly on \,$\Omega$. This complements the results provided in a recent work by the same authors.
\end{abstract}

\maketitle

\section{Introduction}

\quad Throughout this paper, $H(\Omega)$ will stand, as usual, for the vector space of all holomorphic functions on a nonempty open subset \,$\Omega$ \,of the complex plane \,$\C$. It is well known that, if it is endowed with the topology of the uniform convergence on each compact subset of \,$\Omega$ \,(also called compact convergence), then \,$H(\Omega)$ \,becomes a Fr\'echet space, that is, a complete metrizable locally convex topological vector space. Now, if the space \,$H(\Omega)^{\N}$ \,of all sequences of holomorphic functions on \,$\Omega$ \,is equipped with the product topology, then it becomes a Fr\'echet space as well. In the recent work \cite{bernalcalderonlopezprado},
three modes of convergence of sequences \,$(f_n)\in H(\Omega)^{\N}$ \,are considered, namely, pointwise convergence, compact convergence, and uniform convergence on the whole set \,$\Omega$.
Clearly, uniform convergence implies compact convergence, which in turn implies pointwise convergence, both reverse implications being false.

\vskip 3pt

By selecting, for the sake of normalization, the zero as the limit function, it is analyzed in \cite{bernalcalderonlopezprado} the linear size --in an appropriate sense-- of families of sequences \,$\mathbf{f}=(f_n)\in H(\Omega)^{\N}$ \,tending to \,$0$ \,in a given mode of convergence but not in another one. With this aim, the following three sets of sequences in \,$H(\Omega)$ \,are defined in \cite{bernalcalderonlopezprado}:
\begin{itemize}
  \item $\mathcal{S}_p = \left\{\mathbf{f}=(f_n)\in H(\Omega)^{\N} : \, \lim_{n\to\infty}f_n=0 \, \text{ pointwisely}\right\}$.
  \item $\mathcal{S}_{uc} = \left\{\mathbf{f}=(f_n)\in H(\Omega)^{\N} : \, \lim_{n\to\infty}f_n=0 \, \text{ uniformly on compacta}\right\}$.
  \item $\mathcal{S}_{u} = \left\{\mathbf{f}=(f_n) \in H(\Omega)^{\N} : \, \lim_{n\to\infty}f_n=0 \, \text{ uniformly on $\Omega$}\right\}$.
\end{itemize}
As a special case of the above chain of implications, we have that \,$\mathcal{S}_u \subset \mathcal{S}_{uc} \subset \mathcal{S}_p$. Under our terminology, a number of results dealing with the algebraic and algebraic-topological structure of the sets \,$\mathcal{S}_p \setminus \mathcal{S}_{uc}$ \,and \,$\mathcal{S}_{uc} \setminus \mathcal{S}_{u}$ \,have been obtained
in \cite{bernalcalderonlopezprado}. Grosso modo, they assert the existence of large linear vector subspaces or algebras inside both difference sets. These results will be quoted in Section 2, which will also be devoted to recall a number of concepts and results coming from the theory of lineability, as well as to fix the algebraic and topological structure of the set of sequences \,$H(\Omega)^{\N}$.

\vskip 3pt

In this work we complement the results obtained in \cite{bernalcalderonlopezprado} about the linear size of \,$\mathcal{S}_p \setminus \mathcal{S}_{uc}$ \,and \,$\mathcal{S}_{uc} \setminus \mathcal{S}_u$. Specifically, it will be proved in Section 3 that, under the natural topology and with the sole exception of the zero function, we can find a
\emph{closed infinite dimensional vector subspace} of \,$H(\Omega)^{\N}$ \,inside each of these families, even containing respective prescribed members.

\section{Preliminaries and terminology}\label{background}

\quad In the vector space \,$H(\Omega)^{\N}$ \,of all sequences \,$(f_n)$ \,of holomorphic functions in \,$\Omega$, we define the sum and the multiplication by complex scalars in the usual way, and the product of two elements \,${\bf f}=(f_n)$ \,and \,${\bf g}=(g_n)$ \,by \,${\bf f}\cdot{\bf g}=(f_n\cdot g_n)$. Endowed with these operations, $H(\Omega)^{\N}$ \,becomes a commutative linear algebra. Recall that \,$H(\Omega)$ \,has been endowed with the usual compact open topology and then \,$H(\Omega)$ \,becomes a Fr\'echet space (see, e.g., \cite[Section 1.45]{RudinFunctionalAnalysis}). Moreover, thanks to Runge's theorem (see, e.g., \cite[Theorem 13.9]{Rudin}), it is separable (see \cite[pp.~370 and 373]{Kothe}). We equip the space \,$H(\Omega)^{\N}$ \,with its natural topology, that is, the product topology. Then \,$H(\Omega)^{\N}$ \,is also a separable Fr\'echet space.

\vskip 3pt

Next, let us introduce a number of concepts taken from the modern theory of li\-nea\-bi\-li\-ty, for whose background the reader is referred to \cite{ABPS}. A subset \,$A$ \,of a vector space \,$X$ \,is called {\it lineable} whenever there is an infinite dimensional vector subspace of \,$X$ \,that is contained, except for zero, in \,$A$; and \,$A$ \,is said to be {\it algebrable} if it is contained in some linear algebra and there is an infinitely generated algebra contained, except for zero, in \,$A$. In addition, if \,$A$ \,is contained in some commutative linear algebra and \,$\alpha$ \,is a cardinal number, then \,$A$ \,is called {\it strongly $\alpha$-algebrable} if there exists an $\alpha$-generated free algebra \,$M$ \,with \,$M\setminus \{0\} \subset A$.

\vskip 3pt

Now, assume that \,$X$ \,is a topological vector space and \,$A \subset X$. Then \,$A$ \,is called {\it dense-lineable} in \,$X$ \,if there is a dense vector subspace
\,$M\subset X$ \,such that \,$M \setminus\{0\}\subset A$. If there is a closed infinite dimensional vector subspace \,$M\subset X$ \,such that \,$M\setminus\{0\} \subset A$, then \,$A$ \,is called {\it spaceable}.

\vskip 3pt

The more accurate notions of pointwise lineable set and infinitely pointwise lineable set have been recently introduced in \cite{PellegrinoRaposo} and \cite{CalderonGerlachPrado}, respectively. Given a cardinal number \,$\alpha$, a subset \,$A$ \,of a vector space \,$X$ \,is called {\it pointwise $\alpha$-lineable} if for every \,$x \in A$ \,there exists a vector subspace \,$M_x\subset X$ \,such that \,$\dim(M_x)=\alpha$ \,and \,$x\in M_x \subset A \cup \{0\}$. If \,$X$ \,is a topological vector space, then \,$A$ \,is said to be {\it infinitely pointwise $\alpha$-dense-lineable} in \,$X$ \,if for every \,$x \in A$ \,there is a family \,$\left\{W_k : k \in \N \right\}$ \,such that each \,$W_k$ \,is a dense subspace of \,$X$, $\dim(W_k)=\alpha$, $x\in W_k\subset A \cup \{0\}$, and \,$W_k \cap W_n = {\rm span}\{x\}$ \,whenever \,$k \neq n$. Of course, this property is stronger than both dense-lineability and pointwise $\alpha$-lineability. Finally, a subset \,$A$ \,of a topological vector space \,$X$ \,is called {\it pointwise spaceable} if for every \,$x \in A$ \,there exists a closed infinite dimensional vector subspace \,$M_x \subset X$ \,such that \,$x \in M_x \subset A \cup \{0\}$. Clearly, this property is stronger than mere spaceability.

\vskip 3pt

The following theorem gathers all findings given in \cite{bernalcalderonlopezprado} concerning the lineability of \,$\mathcal{S}_p \setminus \mathcal{S}_{uc}$
\,and \,$\mathcal{S}_{uc} \setminus \mathcal{S}_{u}$. The symbol \,$\mathfrak{c}$ \,will represent, as usual, the cardinality of the continuum.

\begin{theorem}\label{Thm del otro paper}
Assume that \,$\Omega$ \,is a nonempty open subset of \,$\C$. In \,$H(\Omega)^{\N}$, we consider the corresponding families of null sequences \,$\mathcal{S}_{p}$, $\mathcal{S}_{uc}$, and \,$\mathcal{S}_{u}$. Then the following holds:
\begin{enumerate}
  \item[\rm (a)] The sets \,$\mathcal{S}_p \setminus \mathcal{S}_{uc}$ \,and \,$\mathcal{S}_{uc} \setminus \mathcal{S}_u$ \,are strongly $\mathfrak{c}$-algebrable.
  \item[\rm (b)] The sets \,$\mathcal{S}_p \setminus \mathcal{S}_{uc}$ \,and \,$\mathcal{S}_{uc} \setminus \mathcal{S}_u$ \,are infinitely pointwise $\mathfrak{c}$-dense-lineable.
  \item[\rm (c)] For each one of the sets \,$\mathcal{E} \in \{ \mathcal{S}_p \setminus \mathcal{S}_{uc}, \mathcal{S}_{uc}\setminus\mathcal{S}_u \}$ \,there exists an infinite dimensional Banach space \,$X \subset H(\Omega)^{\N}$ \,satisfying that \,$X \setminus \{0\} \subset \mathcal{E}$ \,and the norm topology on \,$X$ \,is stronger than the one inherited from \,$H(\Omega)^{\N}$.
\end{enumerate}
\end{theorem}

Observe that, in spite of (c), the {\it spaceability} of our two special families is missed, because nothing indicates that the subspace \,$X$ \,is closed for the product topology on $H(\Omega)^{\N}$. This paper is conceived to fill in this gap.

\section{Spaceability of $\mathcal{S}_p \setminus \mathcal{S}_{uc}$ and $\mathcal{S}_{uc} \setminus \mathcal{S}_u$}

\quad Given a function \,$\varphi:\Omega\to\C$ \,and a subset \,$A\subset\Omega$, then \,$\varphi|_A$ \,will stand, as usual, for the restriction of \,$\varphi$ \,to \,$A$. In addition, $\Omega = \bigcup_{i \in N(\Omega)}\Omega_i$ will be the decomposition of \,$\Omega$ \,into its connected components, where $N(\Omega)$ is either \,$\N$ \,or the set \,$\{1,2,\ldots,p\}$ \,for some \,$p\in\N$. When \,$N(\Omega)=\{1\}$ \,or, which is the same, when \,$\Omega$ \,is connected, then \,$\Omega$ \,is said to be a {\it domain.} We say that a sequence \,$(z_n) \subset \Omega$ \,{\it tends to the boundary} of \,$\Omega$ \,if for each compact subset \,$K \subset \Omega$ \,there is \,$n_0 \in \N$ \,such that \,$z_n \in \Omega \setminus K$ \,for all \,$n\geq n_0$.

\vskip 3pt

Prior to establish the existence of large closed subspaces inside \,$\mathcal{S}_{p} \setminus \mathcal{S}_{uc}$ \,and \,$\mathcal{S}_{uc} \setminus \mathcal{S}_u$, the following two auxiliary lemmas are needed.

\begin{lemma}\label{Lemma H(Omega) as a product}
Under the notation above, the following holds:
\begin{enumerate}
  \item[\rm (a)] If \,$(f_n) \in \mathcal{S}_p \setminus \mathcal{S}_{uc}$, then there exist \,$\alpha>0$, $i\in N(\Omega)$, a sequence \,$(z_n)$ \,of pairwise different points in \,$\Omega_i$, and a subsequence \,$(F_n)$ \,of \,$(f_n)$ \,such that \,$(z_n)$ \,tends to a point \,$z_0 \in \Omega_i$ \,and \,$|F_n(z_n)| \geq \alpha$ \,for all $n\in\N$.

  \item[\rm (b)] If \,$(f_n) \in \mathcal{S}_{uc} \setminus \mathcal{S}_{u}$, then there exist \,$\alpha>0$, a sequence \,$(z_n)$ \,of pairwise different points in \,$\Omega$, and a subsequence \,$(F_n)$ \,of \,$(f_n)$ \,with \,$|F_n(z_n)| \geq \alpha$ \,for all \,$n\in\N$ \,satisfying one of the following properties:
      \begin{itemize}
        \item[\rm (i)] There is \,$i\in N(\Omega)$ \,such that \,$(z_n)$ \,is contained in \,$\Omega_i$ \,and tends to the boundary of \,$\Omega_i$.
        \item[\rm (ii)] $N(\Omega)=\N$ \,and there is a sequence \,$(i(n))$ \,of pairwise different positive integers such that \,$z_n \in \Omega_{i(n)}$ \,for all \,$n\in\N$.
      \end{itemize}
\end{enumerate}
\end{lemma}

\begin{proof}

\noindent (a) Let \,$(f_n) \in \mathcal{S}_p \setminus \mathcal{S}_{uc}$. Then there is a compact set \,$K \subset \Omega$ \,such that \,$(f_n)$ \,converges to zero pointwisely but not uniformly on \,$K$. Hence, there are an $\alpha > 0$, an increasing sequence of positive integers
\,$m(1)<m(2)< \cdots < m(n) < \cdots$ \,and a countable set \,$\{w_n : n \in \N\} \subset K$ \,such that \,$|f_{m(n)}(w_n)| > \alpha$ \,for all \,$n\in\N$. Now, given \,$n\in\N$, the point \,$w_n$ \,cannot appear infinitely many times in the sequence \,$\{w_1,w_2,w_3, \dots \}$ (say \,$w_n = w_{s(j)}$, with \,$s(1) < s(2) < s(3) < \cdots$) because, if this were the case, then we would have \,$|f_{m(s(j))} (w_n)| \geq \alpha$ \,for all \,$j\in\N$ \,and, thus, $\lim_{k\to \infty} f_k(w_n) \neq 0$, which would contradict the pointwise convergence to \,$0$. By discarding the repeated values of the \,$w_n$'s \,and after renaming, we obtain that the \,$w_n$'s \,are pairwise different points of \,$K$ \,and \,$|f_{m(n)}(w_n)| \geq \alpha$ \,for all \,$n \in \N$. Since \,$\Omega = \bigcup_{i\in N(\Omega)}\Omega_i$, the connected components form an open covering of \,$K$. By the compactness of \,$K$, there are finitely many components of \,$\Omega$ \,covering \,$K$. Consequently, at least one of them, say \,$\Omega_i$, must contain infinitely many points of the set \,$\{w_n : \, n \in \N\}$. If \,$w_{n(k)}$ ($k=1,2,3\dots$) are such pairwise different points, then by defining
\begin{equation*}
  F_k := f_{m(n(k))} \hbox{\, and \,} \ z_k := w_{n(k)} \quad (k\in \N),
\end{equation*}
we reach the conclusion.

\vskip 3pt

\noindent (b) Let \,$(f_n) \in \mathcal{S}_{uc} \setminus \mathcal{S}_{u}$. Then \,$f_n\to 0$ compactly but not uniformly on \,$\Omega$. Therefore, there are an \,$\alpha>0$, a sequence of positive integers \,$m(1)<m(2)<\cdots<m(n)<\cdots$ \,and a countable set \,$\{w_n : n \in \N\} \subset \Omega$ \,such that \,$|f_{m(n)}(w_n)| \geq \alpha$ \,for all \,$n\in \N$. Now, given \,$n\in\N$, the point \,$w_n$ \,cannot appear infinitely many times in the sequence \,$\{w_1,w_2,w_3\dots \}$ (say \,$w_n = w_{s(j)}$ \,with \,$s(1) < s(2) < s(3) < \cdots$) because, if this were the case, then we would have \,$|f_{m(s(j))} (w_n)| \geq \alpha$ \,for all \,$j\in\N$ \,and, thus, $\lim_{k\to \infty} f_k(w_n) \neq 0$, which would contradict the compact convergence to \,$0$. By discarding the repeated values of the \,$w_n$'s \,and after renaming, we obtain that the \,$w_n$'s \,are pairwise different points of \,$\Omega$ \,and \,$|f_{m(n)}(w_n)| \geq \alpha$ \,for all \,$n\in\N$. At this point, two cases are possible:
\begin{itemize}
  \item There exists a connected component \,$\Omega_i$ \,containing a subsequence \,$(z_k):= (w_{n(k)})$ \,of \,$(w_n)$. Then we define the subsequence $(F_k):=(f_{m(n(k))})$. Of course, $|F_k(z_k)|\geq \alpha$ for all $k\in\N$. From the assumption \,$(f_n)\in \mathcal{S}_{uc}$, it follows that \,$F_k \to 0$ \,compactly on \,$\Omega_i$ \,and, consequently, each compact subset of \,$\Omega_i$ \,can contain only finitely many \,$z_k$'s. This tells us that \,$(z_k)$ \,tends to the boundary of \,$\Omega_i$, so we obtain (i).

  \item There do not exist a component of \,$\Omega$ \,containing infinitely many points \,$w_n$, that is, $\Omega_i \cap \{w_n : \, n \in \N \}$ \,is finite for all \,$i \in N(\Omega)$. Since \,$\{w_n : \, n \in \N \}$ \,is infinite and contained in \,$\Omega$, we deduce that \,$N(\Omega)=\N$. Now, we proceed by induction. Let \,$n(1):=1$ \,and let \,$i(1)$ \,be the first \,$n\in\N$ \,with \,$w_{n(1)} \in \Omega_n$. If \,$n(1),i(1), \dots ,n(k),i(k)$ \,have been constructed with \,$n(1) < n(2) < \cdots < n(k)$, $i(1) < i(2) < \cdots < i(k)$ \,and \,$w_{n(j)} \in \Omega_{i(j)}$ \,for all \,$j\in \{1,\dots,k\}$, then we define \,$n(k+1)$ \, as the first \,$n\in\N$ \,such that \,$w_n \not\in \Omega_1 \cup \cdots \cup \Omega_{i(k)}$ \,and \,$i(k+1)$ \,as the first \,$n\in\N$ \,with \,$w_{n(k+1)} \in \Omega_n$. Then \,$n(k)<n(k+1)$, $i(k) < i(k+1)$, and \,$w_{n(k+1)} \in \Omega_{i(k+1)}$. Thus, we have constructed two sequences \,$(n(k))$ \,and \,$(i(k))$ \,of natural numbers. Finally, if we take \,$z_k := w_{n(k)}$ \,and \,$F_k := f_{m(n(k))}$ \,for each \,$k \in \N$, then we obtain (ii), because we also have \,$|F_k(z_k)| \geq \alpha$ \,for all \,$k\in\N$.
\end{itemize}
This concludes the proof.
\end{proof}

\vskip 3pt

As usual, the symbol \,$\chi_S$ \,will denote the characteristic function of a set \,$S \subset \Omega$, that is, $\chi_S(z)=1$ \,if \,$z\in S$ \,and $\chi_S(z)=0$ \,if \,$z \not\in S$. Observe that if \,$S$ \,is a union of several connected components of \,$\Omega$, then \,$\chi_S \in H(\Omega)$. Moreover, $\partial A$ \,will denote the boundary of a set \,$A\subset \C$.

\begin{lemma} \label{Lemma para espaciabilidad de Sp-Suc y Suc-Su}
Let \,$X$ be a vector subspace of \,$H(\Omega)$. For each \,${\bf f} = (f_n) \in H(\Omega)^{\N}$, consider the set
\begin{equation*}
  M({\bf f}) := \big\{ (f_n \cdot \varphi) : \varphi \in X \big\}.
\end{equation*}
Then the following properties hold:
\begin{enumerate}
  \item[\rm (a)] $M({\bf f})$ \,is a vector subspace of \,$H(\Omega)^{\N}$.
  \item[\rm (b)] If \,${\bf f}\in\mathcal{S}_p$ \,then \,$M({\bf f})\subset \mathcal{S}_p$.
  \item[\rm (c)] If \,${\bf f}\in\mathcal{S}_{uc}$ \,then \,$M({\bf f})\subset \mathcal{S}_{uc}$.
  \item[\rm (d)] Assume that at least one of the two following conditions is fulfilled:
      \begin{itemize}
        \item[\rm (i)] There exist \,$m\in\N$, a connected component \,$\Omega_i$ \,of \,$\Omega$, and a closed infinite dimensional subspace \,$X_i \subset H(\Omega_i)$ \,not containing the constant function \,$1$ \,such that \,$f_m|_{\Omega_i} \neq 0$ \,and
            \begin{equation*}
              \phantom{aaaaaaa} X = \big\{\lambda + \psi \in H(\Omega) :
              \lambda \in \C, \, \psi|_{\Omega_i} \in X_i, \hbox{ and }
              \psi|_{\Omega_j}=0 \hbox{ for all }  j \in N(\Omega) \setminus \{i\} \big\}.
            \end{equation*}
        \item[\rm (ii)] There is an infinite set \,$\left\{S_k:k\in\N\right\}$ \,of pairwise disjoint nonempty open subsets of \,$\Omega$ \,such that \,$\Omega = \bigcup_{k\in \N} S_k$, for each \,$k\in\N$ \,there is $m(k)\in\N$ satisfying \,$f_{m(k)}|_{S_k} \neq 0$, and
            \begin{equation*}
              X = \left\{\sum_{k=1}^{\infty} c_k\cdot \chi_{S_k} : c_k\in\C \text{ for each } k\in\N \right\}.
            \end{equation*}
      \end{itemize}
      Then \,$M({\bf f})$ \,is closed, infinite dimensional, and \,${\bf f} \in M({\bf f})$.
\end{enumerate}
\end{lemma}

\begin{proof}
(a) This is plain because \,$X$ \,is a vector subspace of \,$H(\Omega)$.

\vskip 3pt

\noindent (b) Given \,$z\in\Omega$, since \,$\lim_{n\to\infty}f_n(z)=0$, then we also have that \,$\lim_{n\to\infty}f_n(z)\varphi(z)=0$ \,for all \,$\varphi\in H(\Omega)$ \,and so for all \,$\varphi \in X$. Thus, the result follows.

\vskip 3pt

\noindent (c) Given a compact set \,$K \subset \Omega$, the continuity of any fixed \,$\varphi\in X$ \,implies that \,$\sup_{z\in K}|\varphi(z)|<\infty$. If \,${\bf f}\in\mathcal{S}_{uc}$, then
\begin{equation*}
  0 \leq \sup_{z\in K} |f_n(z) \cdot \varphi(z)|
  \leq \sup_{z\in K}|\varphi(z)| \cdot \sup_{z\in K}|f_n(z)| \longrightarrow 0  \hbox{\, as \,} n\to\infty
\end{equation*}
and the result follows.

\vskip 3pt

\noindent (d) It is plain that \,$X$ \,is a vector subspace in both cases (i) and (ii). The fact that \,${\bf f} \in M({\bf f})$ \,is derived by choosing \,$\lambda=1$ \,and \,$\psi=0$ \,if (i) holds, and \,$c_k=1$ \,for all \,$k\in\N$ \,if (ii) holds.

\vskip 3pt

Concerning the dimension of \,$M({\bf f})$, let us start with the case (i). Since the function \,$f_m$ \,is not identically zero on the domain \,$\Omega_i$, it follows from the
Identity Principle that if \,$\varphi \in H(\Omega_i)$ \,and \,$f_m \cdot \varphi=0$ \,on \,$\Omega_i$, then \,$\varphi=0$. It follows that if \,$\left\{\varphi_k:k\in\N\right\}$ \,is an infinite set of linearly independent functions from \,$X_i$, then \,$f_m\cdot\varphi_1, f_m\cdot \varphi_2, f_m\cdot\varphi_3,\ldots$ \,are also linearly independent in $H(\Omega_i)$. Each \,$\varphi_k$ \,is now extended to the whole \,$\Omega$ \,by defining it as \,$0$ \,on \,$\Omega \setminus \Omega_i$. Then the sequences \,$(f_n\cdot\varphi_1),(f_n\cdot\varphi_2),(f_n\cdot\varphi_3),\ldots$ \,are linearly independent members of \,$M({\bf f})$, which proves that this space has infinite dimension.

\vskip 3pt

For the case (ii), the sequences \,$(f_n\cdot \chi_{S_1}),(f_n\cdot \chi_{S_2}),(f_n\cdot \chi_{S_3}),\ldots$ \,are linearly independent. Indeed, suppose that \,$k\in\N$, $\lambda_1,\ldots,\lambda_k\in\C$, and
\begin{equation*}
  \lambda_1(f_n\cdot \chi_{S_1}) + \cdots + \lambda_k(f_n\cdot \chi_{S_k}) = 0.
\end{equation*}
The disjointness of the sets \,$S_k$ \,implies that \,$\lambda_1f_n=0$ \,on $S_1$ \,for every \,$n\in\N$. By (ii), there exists \,$m(1)\in\N$ \,such that \,$f_{m(1)}\neq 0$ \,on \,$S_1$, so \,$\lambda_1=0$. If the argument is repeated, we obtain that \,$\lambda_1=\cdots=\lambda_k=0$. This shows that \,${\rm dim}(M({\bf f}))=\infty$ also in case (ii).

\vskip 3pt

Let us prove that \,$M({\bf f})$ \,is closed in \,$H(\Omega)$ under any of assumptions (i) or (ii).
With this aim, assume that \,$({\bf g}_k)$ \,is a sequence of members of \,$M({\bf f})$ \,such that \,${\bf g}_k \to {\bf h} = (h_n)\in H(\Omega)^{\N}$ \,as \,$k\to\infty$ \,in the product topology. We have to show that \,${\bf h} \in M({\bf f})$; that is, we have to prove the existence of a function \,$\Phi\in X$ \,such that \,$h_n = f_n \cdot \Phi$ \,for all \,$n\in\N$. For each \,$k\in\N$, the sequence \,${\bf g}_k$ \,has the form \,${\bf g}_k = (f_n \cdot \varphi_k)$ \,for some \,$\varphi_k \in X$. Then, for each coordinate \,$n$, we have that
\begin{equation}\label{Eq:fn phik tiende a hn}
  f_n \cdot \varphi_k \longrightarrow h_n \hbox{\, compactly on \,} \Omega \hbox{\, as \,} k\to\infty.
\end{equation}

Assume that the property (i) holds. In this case there exist sequences \,$(\lambda_k)\subset\C$ \,and \,$(\psi_k)\subset H(\Omega)$ \,such that \,$\varphi_k=\lambda_k+\psi_k$,
$\psi|_{\Omega_i} \in X_i$, and \,$\psi_k|_{\Omega_j}=0$ \,for all $j \neq i$. Observe that each function \,$\varphi_k|_{\Omega_i}$ \,belongs to the subspace \,$Y := \langle 1 \rangle \oplus X_i$, the algebraic direct sum (because \,$1 \not\in X_i$) of the linear span of \,$1$ (that is, the space of constants) and \,$X_i$. Since \,$\langle 1 \rangle$ \,is finite dimensional, the subspace \,$Y$ \,is closed in \,$H(\Omega_i)$ (see, e.g., \cite[Theorem 1.42]{RudinFunctionalAnalysis}). In addition, \,$Y$ \,is also the {\it topological} direct sum of \,$\langle 1 \rangle$ \,and \,$X_i$, that is, the projections
\begin{equation*}
  \pi_1 : \lambda + \psi \in \langle 1 \rangle \oplus X_i \longmapsto \lambda \in \C
  \quad \hbox{and} \quad
  \pi_2 : \lambda + \psi \in \langle 1 \rangle \oplus X_i \longmapsto \psi \in X_i
\end{equation*}
are continuous (see \cite[p.~22]{SchaeferWolff}).

\vskip 3pt

Let \,$z_0\in\Omega_i$ \,and let \,$\mu\in\N\cup\{0\}$ \,be the order of \,$z_0$ \,as a zero of the function \,$f_m$ \,that appears in the property (i). Then there exists a function \,$F\in H(\Omega_i)$ \,such that \,$F(z_0)\neq 0$ \,and $f_m(z)=(z-z_0)^\mu \cdot F(z)$ \,for all \,$z\in \Omega_i$. In addition, there exist two constants $R>0$ and $\alpha>0$ such that the closed disc \,$D=\{z\in\C:|z-z_0|\leq R\}$ \,is contained in \,$\Omega_i$ \,and \,$|F(z)|\geq \alpha$ \,for all $z\in D$. Since the circle \,$\partial D$ \,is a compact subset of \,$\Omega_i$, the property \eqref{Eq:fn phik tiende a hn} implies the following:
\begin{align*}
  \sup_{z\in\partial D} \left|\varphi_k(z) - \frac{h_m(z)}{f_m(z)} \right|
  &= \sup_{z\in\partial D}\frac{1}{|z-z_0|^\mu \cdot |F(z)|} \cdot |f_m(z)\varphi_k(z) - h_m(z)| \\
  &\leq \frac{1}{R^\mu \cdot \alpha} \cdot \sup_{z\in\partial D}|f_m(z) \varphi_k(z) - h_m(z)|  \longrightarrow 0 \hbox{ as } k \to \infty.
\end{align*}
Consequently, $\sup_{z\in\partial D}\left|\varphi_k(z)-\frac{h_m(z)}{f_m(z)}\right|\to 0$ \,as \,$k\to\infty$. Given \,$\varepsilon>0$, we derive the existence of \,$k_0\in\N$ \,such that \,$\sup_{z\in\partial D}\left|\varphi_k(z)-\frac{h_m(z)}{f_m(z)}\right|< \varepsilon/2$ \,for all
\,$k\geq k_0$. It follows from the triangle inequality that
\,$\sup_{z\in\partial D}|\varphi_k(z)-\varphi_l(z)|<\varepsilon$ \,for all \,$k,l\geq k_0$. Now, we invoke the Maximum Modulus Principle to deduce that \,$\sup_{z \in D}\left|\varphi_k(z)-\varphi_l(z)\right|<\varepsilon$ \,for all \,$k,l\geq k_0$. Since every compact subset of \,$\Omega_i$ \,can be covered by a finite amount of discs contained in \,$\Omega_i$, we deduce that \,$(\varphi_k)$ is a Cauchy sequence in the Fr\'echet space $H(\Omega_i)$. Consequently, there is \,$\Phi_i \in H(\Omega_i)$ \,such that
\begin{equation} \label{varphik tiende a Phii}
  \varphi_k \longrightarrow \Phi_i  \hbox{\, compactly on \,} \Omega_i  \hbox{\, as \,} k\to\infty.
\end{equation}
Since \,$Y$ \,is closed, we have \,$\Phi_i \in Y$. Furthermore, since compact convergence implies
pointwise convergence and the pointwise limit is unique, we get \,$\Phi_i(z)=\frac{h_m(z)}{f_m(z)}$ \,for all \,$z\in\Omega_i \setminus Z$, where \,$Z:=\left\{z\in\Omega_i:f_m(z)=0\right\}$ \,is a discrete subset of \,$\Omega_i$ \,due to the Identity Principle. This entails that \,$h_m(z)=f_m(z)\Phi_i(z)$ \,for every \,$z$ \,in \,$\Omega \setminus Z$, which is a dense subset of \,$\Omega_i$. The continuity of both functions \,$h_m$ \,and \,$f_m \cdot \Phi_i$ \,implies \,$h_m=f_m\cdot\Phi_i$ \,on the whole $\Omega_i$.

\vskip 3pt

On the one hand, since \,$\Phi_i\in Y$, one derives the existence of \,$\lambda\in\C$ \,and \,$\Psi_i \in X_i$ \,such that \,$\Phi_i=\lambda+\Psi_i$. On the other hand, the continuity of \,$\pi_1$ \,and \,$\pi_2$ \,and the fact \eqref{varphik tiende a Phii} imply that \,$\psi_k \to \Psi_i$ \,compactly
on \,$\Omega_i$ \,and \,$\lambda_k \to \lambda$ \,as \,$k\to\infty$. Define the function \,$\Phi \in H(\Omega)$ \,as follows:
\begin{equation*}
  \Phi(z) = \begin{cases}
              \Phi_i(z) = \lambda + \Psi_i(z) & \mbox{if }  z\in \Omega_i \\
              \lambda                         & \mbox{otherwise.}
            \end{cases}
\end{equation*}
Clearly, $\Phi\in X$. Moreover, \eqref{Eq:fn phik tiende a hn} together with the uniqueness of the limit tell us that for each \,$n\in\N$ \,we have \,$h_n=f_n \cdot (\lambda + \Psi)$ \,on \,$\Omega_i$, while \,$h_n = \lambda \cdot f_n$ \,on the remaining components \,$\Omega_j$ \,with \,$j\neq i$. In other words, $h_n = f_n \cdot \Phi$ \,on \,$\Omega$ \,for all \,$n \in \N$. This indicates that \,${\bf h} \in M({\bf f})$, as required.

\vskip 3pt

Finally, suppose that the property (ii) holds. In this case, each function \,$\varphi_k$ \,can be written as \,$\varphi_k = \sum_{\nu=1}^{\infty} c_{k,\nu} \cdot \chi_{S_\nu}$, where \,$c_{k,\nu}\in\C$ \,for all $k,\nu\in \N$. Now, \eqref{Eq:fn phik tiende a hn} yields that, for all \,$n,\nu\in \N$, we have that \,$f_n \cdot \varphi_k \to h_n$ \,uniformly on the compact subsets of \,$S_{\nu}$ \,as \,$k\to \infty$. Fix \,$\nu \in \N$. By (ii), there exist an index \,$m(\nu)\in \N$ \,and a point \,$z_0 \in S_{\nu}$ \,such that \,$f_{m(\nu)}(z_0) \neq 0$. Since uniform convergence implies pointwise convergence and \,$\varphi_k|_{S_\nu} = c_{k,\nu}$, we obtain that \,$f_{m(\nu)}(z_0) \cdot c_{k,\nu} \to h_{m(\nu)}(z_0)$ \,as $k\to\infty$, so
\,$c_{k,\nu}\to \frac{h_{m(\nu)}(z_0)}{f_{m(\nu)}(z_0)}$ \,as $k\to\infty$. Then the value
\begin{equation*}
  d_\nu := \frac{h_{m(\nu)}(z_0)}{f_{m(\nu)}(z_0)} = \lim_{k\to\infty}c_{k,\nu} \in \C
\end{equation*}
is independent of \,$z_0$. Thus, for every pair \,$n,\nu \in \N$ \,we have, on the one hand, that \,$f_n \cdot c_{k,\nu} \to h_n$ \,compactly on \,$S_\nu$ \,as \,$k\to\infty$ \,and, on the other hand, that \,$f_n \cdot c_{k,\nu} \to f_n \cdot d_\nu$ \,compactly on \,$S_\nu$ \,as \,$k\to\infty$. The uniqueness of the limit implies that \,$h_n=f_n \cdot d_\nu$ \,on \,$S_\nu$ \,for all \,$n,\nu\in\N$. Recall that \,$\Omega = \bigcup_{\nu=1}^{\infty}S_\nu$ \,and that this is a disjoint union of open sets. Consequently, if we define \,$\Phi := \sum_{\nu=1}^{\infty} d_\nu \cdot \chi_{S_\nu} \in X$, then we get, as in the preceding case, that \,$h_n = f_n \cdot \Phi$ \,on \,$\Omega$ \,for all \,$n\in\N$, which concludes the proof.
\end{proof}

\vskip 3pt

\begin{theorem} \label{Thm-Sp-Suc spaceable}
The set \,$\mathcal{S}_{p} \setminus \mathcal{S}_{uc}$ \,is pointwise spaceable in \,$H(\Omega)^{\N}$.
\end{theorem}

\begin{proof}
First of all, note that the set $\mathcal{S}_p \setminus \mathcal{S}_{uc}$ is not empty by Theorem \ref{Thm del otro paper}. Let \,${\bf f} = (f_n) \in \mathcal{S}_p \setminus \mathcal{S}_{uc}$. By Lemma \ref{Lemma H(Omega) as a product}, there exist $\alpha > 0$, $i \in N(\Omega)$, a sequence \,$(z_n)$ \,of pairwise different points in \,$\Omega_i$, and a subsequence \,$(F_n)$ \,of \,$(f_n)$ \,such that \,$(z_n)$ \,tends to a point \,$z_0 \in \Omega_i$ \,and \,$|F_n(z_n)|\geq \alpha$ \,for all $n\in\N$.

\vskip 3pt

Select and fix any point \,$a \in \Omega_i \setminus \{z_0,z_1,z_2, \dots \}$. The continuity of the evaluation functional \,$\psi \in H(\Omega) \mapsto \psi(a)\in\C$ \,implies that the set
\begin{equation*}
  X_i := \{\psi \in H(\Omega_i) : \psi(a)=0\}
\end{equation*}
is a closed subspace of \,$H(\Omega_i)$. Moreover, $X_i$ \,is  infinite dimensional because it contains the functions given by \,$\psi_n(z)=(z-a)^n$ \,for all \,$n\in\N$. Consider the space
\begin{equation*}
  X :=  \big\{\lambda + \psi \in H(\Omega) : \lambda \in \C, \, \psi|_{\Omega_i} \in X_i \,
  \hbox{ and } \, \psi|_{\Omega_j}=0 \, \hbox{ for all } \, j\in N(\Omega)\setminus \{i\} \big\}.
\end{equation*}
Since condition (i) in Lemma \ref{Lemma para espaciabilidad de Sp-Suc y Suc-Su} is fulfilled, we derive that the set
\begin{equation*}
  M({\bf f}) := \big\{ (f_n \cdot \varphi) : \varphi \in X \big\}
\end{equation*}
is a closed infinite dimensional subspace of \,$H(\Omega)^\N$ \,containing \,${\bf f}$ \,and contained in \,$\mathcal{S}_p$. It remains to show that any \,${\bf g} = (g_n) \in M({\bf f}) \setminus \{0\}$ \,does not belong to \,$\mathcal{S}_{uc}$.

\vskip 3pt

With this aim, notice that for such a \,${\bf g}$ \,there is \,$\varphi \in X \setminus \{0\}$ \,such that \,$g_n = f_n \cdot \varphi$ \,for all \,$n\in\N$. Then there are \,$\lambda \in \C$ \,and \,$\psi \in X_i$ \,such that \,$\varphi = \lambda + \psi$. Assume, by way of contradiction, that \,$\varphi=0$ \,on $\Omega_i$. Since \,$\psi(a)=0$, we would have that \,$\lambda=\varphi(a)=0$, so \,$\varphi|_{\Omega_j}=\psi|_{\Omega_j}=0$ \,for all $j\neq i$. Thus, \,$\varphi=0$ \,on \,$\Omega$, which is not true. Therefore, $\varphi$ \,cannot be identically zero on \,$\Omega$. By the Identity Principle, there is \,$r>0$ \,such that the disc $K:=\left\{z\in\C:|z-z_0|\leq r\right\}$ \,is contained in \,$\Omega_i$ \,and \,$\varphi(z)\neq 0$ \,for all \,$z\in \partial K$. In particular, there exists $\beta \in (0,+\infty)$ such that $|\varphi(z)| \geq \beta$ for all $z \in \partial K$. Since \,$z_n \to z_0$, there is \,$n_0 \in \N$ \,such that \,$z_n \in K$ \,for all \,$n \geq n_0$. Finally, we invoke the Maximum Modulus Principle to obtain that
\begin{equation*}
  0 < \alpha \leq |F_n(z_n)| \leq \sup_{z\in K} |F_n(z)| = \sup_{z\in \partial K}|F_n(z)|
  \leq {1 \over \beta} \cdot \sup_{z\in \partial K}|F_n(z) \cdot \varphi(z)|
\end{equation*}
for all \,$n \geq n_0$. This proves that \,$(F_n \cdot \varphi)$ \,does not converges to zero uniformly on the compact set \,$\partial K$. Taking into account that \,$(F_n)$ \,is a subsequence of \,$(f_n)$, it follows that \,${\bf g}=(f_n \cdot \varphi) \not\in\mathcal{S}_{uc}$.
\end{proof}

\vskip 3pt

Before proving the spaceability of the set \,$\mathcal{S}_{uc} \setminus \mathcal{S}_{u}$, we need some preparation about basic sequences in Banach spaces. The symbol \,$\D$ \,will stand for the open unit disc in the complex plane, \,$\T = \partial \D$ \,will denote the unit circle, and \,$L^2(\T)$ \,will represent the Hilbert space \,of all Lebesgue classes of measurable functions \,$f:\T \to \C$ \,with finite quadratic norm:
\begin{equation*}
  \|f\|_2 :=
  \left(\frac{1}{2\pi}\int_0^{2\pi}\left|f\left(e^{i\theta}\right)\right|^2 d\theta\right)^{1/2}
  < \infty.
\end{equation*}
Recall that a sequence \,$(x_j)$ \,in a Banach space is said to be a {\em basic sequence} whenever every vector \,$x$ \,in the closed linear span of \,$(x_j)$ \,can be written as \,$x=\sum_{j=1}^\infty a_j x_j$ \,for a unique scalar sequence \,$(a_j)$. Two basic sequences \,$(x_j)$ \,and \,$(y_j)$ \,are said to be {\it equivalent} if, for every sequence
$(a_j)$ \,of scalars, the series \,$\sum_{j=1}^{\infty}a_j x_j$ \,converges if and only if the series \,$\sum_{j=1}^{\infty}a_j y_j$ \,converges. Observe that \,$(z^j)$ \,is a basic
sequence in \,$L^2(\T)$ \,because \,$\left\{z^j : j\in\mathbb{Z}\right\}$ \,is an orthonormal basis of \,$L^2(\T)$.

\begin{lemma}\label{basic sequence}
Assume that \,$\Omega$ is an open set of \,$\C$ \,with \,$\overline{\D} \subset \Omega$ \,and that
\,$(\Phi_j)\subset H(\Omega)$ \,is a basic sequence in \,$L^2(\T)$ \,equivalent to \,$(z^j)$. Let \,$\left\{h_l := \sum_{j=1}^{J(l)} c_{j,l}\Phi_j : l\in\N\right\}$ \,be a sequence in
\,${\rm span }\left\{\Phi_j:j\in\N\right\}$ that converges in $H(\Omega)$. Then
\begin{equation*}
  \sup_{l \in \N} \sum_{j=1}^{J(l)} |c_{j,l}|^2 < +\infty.
\end{equation*}
\end{lemma}

\begin{proof}
Since \,$\T$ \,is a compact subset of \,$\Omega$, convergence in \,$H(\Omega)$ \,is stronger than convergence in \,$L^2(\T)$-norm. Therefore, the sequence \,$(h_l)$ \,converges in $L^2(\T)$ \,and, as a consequence, $\alpha:=\sup_{l\in\N}\|h_l\|_2<\infty$. Since the basic sequences \,$(\Phi_j)$ \,and \,$(z^j)$ \,are equivalent, there are two constants \,$C_1,C_2\in(0,+\infty)$ \,such that
\begin{equation*}
  C_1 \left\|\sum_{j=1}^J a_j z^j\right\|_2 \leq \left\|\sum_{j=1}^J a_j \Phi_j\right\|_2
  \leq C_2 \left\|\sum_{j=1}^J a_j z^j\right\|_2
\end{equation*}
for all \,$J \in \N$ \,and all scalars \,$a_1,\dots,a_J$ (see, e.g., \cite[p.~170]{Fabian}). Thus, for every \,$l \in \N$ \,we have that
\begin{equation*}
  C_1^2 \sum_{j=1}^{J(l)} |c_{j,l}|^2
  = C_1^2 \left\|\sum_{j=1}^{J(l)}c_{j,l} z^j\right\|_2^2
  \leq \left\|\sum_{j=1}^{J(l)} c_{j,l}\Phi_j\right\|_2^2
  = \|h_l\|_2^2 \leq \alpha^2 .
\end{equation*}
Hence, $\sup_{l \in \N} \sum_{j=1}^{J(l)} |c_{j,l}|^2\leq \frac{\alpha^2}{C_1^2}$. That concludes the proof.
\end{proof}

\vskip 3pt

\begin{theorem} \label{Thm-Suc-Su spaceable}
The set \,$\mathcal{S}_{uc} \setminus \mathcal{S}_{u}$ \,is pointwise spaceable in \,$H(\Omega)^{\N}$.
\end{theorem}

\begin{proof}
The set \,$\mathcal{S}_{uc} \setminus \mathcal{S}_{u}$ \,is nonempty by Theorem \ref{Thm del otro paper}. Let \,${\bf f} = (f_n) \in \mathcal{S}_{uc} \setminus \mathcal{S}_{u}$. According to Lemma \ref{Lemma H(Omega) as a product}, there exist \,$\alpha>0$, an infinite sequence \,$(z_n)$ \,in \,$\Omega$ \,of pairwise different points, and a subsequence \,$(F_n)$ \,of \,$(f_n)$ \,with \,$|F_n(z_n)| \geq \alpha$ \,for all $n\in\N$ satisfying one of the following properties:
\begin{itemize}
  \item[\rm (i)] There is \,$i\in N(\Omega)$ \,such that \,$(z_n)$ \,is contained in \,$\Omega_i$ \,and tends to the boundary of \,$\Omega_i$.
  \item[\rm (ii)] $N(\Omega)=\N$ \,and there is a sequence \,$(i(n))$ \,of pairwise different positive integers such that \,$z_n \in \Omega_{i(n)}$ \,for all \,$n\in\N$.
\end{itemize}

Firstly, we face the case in which (i) happens. Without loss of generality, we can assume that \,$\overline{\D}\subset \Omega_i$. Let \,$(\Omega_i)_* := \Omega_i \cup \{\omega \}$ \,be the one-point compactification of \,$\Omega_i$, so that the added point \,$\omega$ \,represents the whole boundary of \,$\Omega_i$ \,in the extended complex plane \,$\C_\infty=\C\cup\{\infty\}$.

\vskip 3pt

Since \,$\overline{\D}\subset \Omega_i$ \,and \,$(z_n)$ \,tends to the boundary of \,$\Omega_i$, by deleting finitely many points \,$z_n$ \,if necessary, we can assume that \,$|z_n|>1$ \,for all \,$n\in\N$. Define the set
\begin{equation*}
  A := \overline{\D} \cup \{z_n : n \in \N\}.
\end{equation*}
Note that \,$A$ \,is a relatively closed subset of \,$\Omega_i$ \,because the sequence \,$(z_n)$ \,tends to the boundary. In addition, the set \,$(\Omega_i)_* \setminus A$ \,is connected as well as locally connected at \,$\omega$  \,because \,$\overline{\D}$ \,is compact (so it is ``far'' from \,$\omega$ \,and we can suppose that the basic connected neighborhoods of \,$\omega$ \,do not
intersect \,$\overline{\D}$), $\Omega_i \setminus \overline{\D}$ \,is connected and \,$(z_n)$
\,is countable (so deleting \,$\left\{z_n: n \in \N\right\}$\, from \,$\Omega_i \setminus \overline{\D}$ \,makes no influence in connectedness or local connectedness).

\vskip 3pt

Now, let us split \,$(z_n)$ \,into infinitely many pairwise disjoint sequences \,$\left\{\left(a_{k,n}\right)_{n=1}^{\infty}:k\in\N\right\}$ \,and consider,
for every \,$k\in \N$, the function \,$g_k:A\to\C$ \,defined as follows:
\begin{equation*}
  g_k(z) =  \begin{cases}
              z^{k-1} & \mbox{if } z\in \overline{\D},  \\
              n       & \mbox{if } z=a_{k,n} \mbox{ for some } n\in \N, \\
              0       & \mbox{if } z=a_{j,n} \mbox{ for some } j,n\in \N \mbox{ with } j\neq k.
            \end{cases}
\end{equation*}
Observe that \,$g_k$ \,is continuous on \,$A$ \,and holomorphic on its interior \,$A^\circ=\D$. Then the Arakelian approximation theorem (see \cite[pp.~136--144]{Gaier}) guarantees the existence of a
function \,$\varphi_k \in H(\Omega_i)$ \,satisfying
\begin{equation*}
  |\varphi_k(z) - g_k(z)| < \frac{1}{3^k}
\end{equation*}
for all \,$z\in A$. Consequently, one obtains
\begin{equation} \label{Eq Arakelian 1}
  |\varphi_k(z) - z^{k-1}| < \frac{1}{3^k} \hbox{\, for all \,} z\in\overline{\D},
\end{equation}
\begin{equation} \label{Eq Arakelian 2}
  |\varphi_k(a_{k,n}) - n| < 1  \mbox{\, for all \,} n\in\N , \mbox{\, and \,}
\end{equation}
\begin{equation} \label{Eq Arakelian 3}
  |\varphi_k(a_{j,n})| < \frac{1}{3^k} \mbox{\, for all \,} n,j\in\N \mbox{\, with \,} j\neq k.
\end{equation}

\vskip 3pt

For each \,$k\in\N$, let \,$\Phi_k\in H(\Omega_i)$ \,stand for the function defined as \,$\Phi_k(z):= z \cdot \varphi_k(z)$ \,for each \,$z\in\Omega_i$. Let \,$(\psi_k)$ \,be the basic sequence in \,$L^2(\T)$ \,given by \,$\psi_k(z):=z^k$ \,for each \,$z\in\T$ \,and let \,$(\psi_k^*)$ be the sequence of coefficient functionals corresponding to \,$(\psi_k)$. Since \,$\|\psi_k^*\|_2 = 1$ \,for all \,$k\in \N$, we have from \eqref{Eq Arakelian 1} that
\begin{equation*}
  \sum_{k=1}^{\infty} \|\psi_k^*\|_2 \cdot \|\Phi_k - \psi_k\|_2
  = \sum_{k=1}^{\infty}\frac{1}{3^k} < 1.
\end{equation*}
By the Nikolskii basis perturbation theorem (see \cite[p.~46, Theorem 9]{Diestel}), we derive that \,$(\Phi_k)$ \,is a basic sequence in \,$L^2(\T)$.

\vskip 3pt

Let us define the set
\begin{equation} \label{Definition of X}
  X_i := \overline{{\rm span} \, \{\Phi_k : k \in \N\}},
\end{equation}
where the closure is considered with respect to the compact open topology on \,$H(\Omega_i)$. Then \,$X_i$ \,is a closed vector subspace of \,$H(\Omega_i)$. Since \,$(\Phi_k|_{\T})$ \,is a basic sequence in \,$L^2(\T)$, their elements are linearly independent as functions on \,$\T$. This implies that the \,$\Phi_k$'s are also linearly independent as functions on \,$\Omega_i$, so \,$X_i$ \,is infinite dimensional. Moreover, $1\not\in X_i$ \,because all functions from \,$X_i$ \,vanish at \,$0$.

\vskip 3pt

Consider the vector space
\begin{equation*}
  X := \big\{\lambda + \varphi \in H(\Omega) : \lambda \in \C, \, \varphi|_{\Omega_i} \in X_i, \hbox{\, and \,} \varphi|_{\Omega_j} = 0 \, \hbox{ for all } \, j \in N(\Omega) \setminus \{i\} \big\}.
\end{equation*}
Our assumption (i) entails that for at least one \,$m\in\N$ it holds that \,$f_m|_{\Omega_i} \neq 0$. It follows from parts (c) and (d) of Lemma \ref{Lemma para espaciabilidad de Sp-Suc y Suc-Su} that
\begin{equation*}
  M({\bf f}) := \big\{ (f_n \cdot \Phi ): \Phi \in X \big\}
\end{equation*}
is a closed infinite dimensional vector subspace of \,$H(\Omega)^{\N}$ \,such that \,${\bf f} \in M({\bf f}) \subset \mathcal{S}_{uc}$. Consequently, it only remains to prove that \,${\bf g} \not\in \mathcal{S}_u$ \,provided that ${\bf g} \in M({\bf f}) \setminus \{0\}$.

\vskip 3pt

Let \,${\bf g}=(g_n)\in M({\bf f}) \setminus \{0\}$. On the one hand, there is \,$\Phi \in X \setminus \{0\}$ \,such that \,$g_n = f_n \cdot \Phi$ \,for all \,$n\in\N$. On the other hand, since \,$(F_n)$ \,is a subsequence of \,$(f_n)$, it suffices to show that some subsequence of \,$(F_n\cdot\Phi)$ \,does not converge to \,$0$ \,uniformly on \,$\Omega$. By the definition of \,$X$, there are \,$\lambda\in\C$ \,and \,$\varphi\in H(\Omega)$ \,such that \,$\varphi|_{\Omega_i} \in X_i$, $\varphi|_{\Omega_j}=0$ \,for all \,$j\neq i$, and \,$\Phi=\lambda+\varphi$. Since convergence with respect to the compact open topology on \,$H(\Omega_i)$ \,is stronger than convergence in \,$L^2(\T)$, we have that \,$\varphi|_{\T}$ \,belongs to
\begin{equation*}
  \widetilde{X} := \overline{{\rm span}\, \{ \Phi_k: k \in \N\})},
\end{equation*}
where the closure is now considered with respect to the norm-topology on \,$L^2(\T)$. Since \,$(\Phi_k)$ \,is a basic sequence in this Banach space, the function \,$\varphi|_{\Omega_i}$ \,has a unique representation \,$\varphi = \sum_{j=1}^{\infty}b_j \Phi_j$ \,in $L^2(\T)$ (with \,$b_j \in \C$ for all \,$j$). Then there exists an increasing sequence \,$p_1<p_2<p_3<\cdots$ \,in \,$\N$ \,such that
\begin{equation*}
  \Phi(z) = \lambda+\varphi(z) = \lambda + \lim_{l\to\infty}\sum_{j=1}^{p_l} b_j \Phi_j(z)
\end{equation*}
for almost every \,$z\in\T$ (see \cite[Theorem 3.12]{Rudin}).

Two cases are possible:
\begin{itemize}
  \item $b_j=0$ for all $j\in\N$.
  \item There exists \,$k\in\N$ \,with \,$b_k\neq 0$.
\end{itemize}
If \,$b_j=0$ \,for all \,$j\in\N$, then \,$\Phi(z)=\lambda$ \,for almost every \,$z\in\T$. By the Identity Theorem, \,$\Phi=\lambda$ \,on \,$\Omega_i$. Since \,$\varphi|_{\Omega_j}=0$ \,for all \,$j\neq i$, it follows that \,$\Phi=\lambda$ \,on the whole \,$\Omega$. Recall that \,$\Phi\neq 0$, so \,$\lambda\neq 0$. Then it is trivial that \,$(f_n \cdot \Phi)=(\lambda f_n)$ \,does not converge to \,$0$ \,uniformly on \,$\Omega_i$, hence \,${\bf g} \not\in \mathcal{S}_u$.

\vskip 3pt

Let us conclude with the analysis of the second case: there exists \,$k\in\N$ (which will be fixed for this part of the proof) with \,$b_k\neq 0$. Observe that since \,$\Phi-\lambda=\varphi$ \,belongs to \,$X_i$, there is a sequence \,$\left\{h_l:=\sum_{j=1}^{J(l)} c_{j,l}\Phi_j : l\in\N\right\}$ \,in ${\rm span}\{\Phi_j: j\in\N\}$ \,such that \,$h_l\to \Phi-\lambda$ \,as \,$\l\to\infty$ \,compactly on \,$\Omega_i$. Without loss of generality, we can assume that \,$J(l)\geq k$ \,for all \,$l$.  By Lemma \ref{basic sequence}, it holds that
\begin{equation}\label{Eq def de C}
  C:= \sup_{l \in \N} \sum_{j=1}^{J(l)} |c_{j,l}|^2 < +\infty .
\end{equation}
But \,$(h_l)$ \,also converges to \,$\Phi-\lambda$ \,in \,$L^2(\T )$, so the continuity of the projection
\begin{equation*}
  \sum_{j=1}^{\infty} d_j \Phi_j \in \widetilde{X} \longmapsto d_k \in \C
\end{equation*}
yields that \,$\lim_{l\to\infty}c_{k,l}=b_k\neq 0$. In particular, there exists \,$l_0\in\N$ \,such that
\begin{equation} \label{Eq coeffs ck}
  |c_{k,l}| > \frac{|b_k|}{2}  \mbox{\, for all \,} l\geq l_0.
\end{equation}

For each \,$n\in\N$, the singleton \,$\{a_{k,n}\}$ \,is a compact subset of \,$\Omega_i$, so there exists a positive integer \,$l_n\geq l_0$ \,such that
\begin{equation} \label{Eq hl}
  |\lambda + h_{l_n}(a_{k,n}) - \Phi(a_{k,n})| < 1.
\end{equation}
Now, we set \,$G_n := F_{n^*}$, where, for each \,$n\in\N$, the symbol \,$n^*$ \,denotes the unique natural number such that \,$z_{n^*}=a_{k,n}$. Then \,$(G_n \cdot \Phi)$ \,is a subsequence of \,$(F_n\cdot \Phi)$, $|a_{k,n}|>1$, and \,$|G_n(a_{k,n})| \geq \alpha$ \,for all \,$n\in\N$. Let \,$\beta :=\alpha+\alpha |\lambda|$. Using \eqref{Eq Arakelian 2}, \eqref{Eq Arakelian 3}, \eqref{Eq coeffs ck}, \eqref{Eq hl}, and the triangle inequality, we obtain the following inequalities:
\begin{align*}
  |G_n(a_{k,n}) \cdot \Phi(a_{k,n})|
  &> \alpha\left(|h_{l_n}(a_{k,n})|-1-|\lambda|\right) = \alpha |h_{l_n}(a_{k,n})| - \beta \\
  &\geq \alpha |c_{k,l_n} \Phi_k(a_{k,n})|
  - \alpha \cdot \sum_{j=1 \atop j \neq k}^{J(l_n)} |c_{j,l_n} \Phi_j(a_{k,n})| - \beta \\
  &= |a_{k,n}| \left( \alpha |c_{k,l_n} \varphi_k(a_{k,n})|
  - \alpha \cdot \sum_{j=1 \atop j \ne k}^{J(l_n)} |c_{j,l_n} \varphi_j(a_{k,n})| \right) - \beta \\
  &> \frac{\alpha |b_k|}{2}(n-1)
  - \alpha \cdot \sum_{j=1 \atop j \ne k}^{J(l_n)} \frac{|c_{j,l_n}|}{3^j} - \beta.
\end{align*}
Now, the Cauchy-Schwarz inequality and \eqref{Eq def de C} imply the following:
\begin{align*}
  |G_n(a_{k,n}) \cdot \Phi(a_{k,n})|
  &> \frac{\alpha |b_k|}{2}(n-1)
  - \alpha \cdot \left(\sum_{j=1}^{J(l_n)} \left(\frac{1}{3^j}\right)^2 \right)^{1/2} \cdot
  \left(\sum_{j=1}^{J(l_n)} |c_{j,l_n}|^2 \right)^{1/2} - \beta \\
  &> \frac{\alpha |b_k|}{2}(n-1) - \alpha \cdot C^{1/2} - \beta
  \longrightarrow +\infty \hbox{\, as \,} n\to\infty.
\end{align*}
Consequently, $\lim_{n \to \infty} \sup_{z\in\Omega_i} |G_n(z) \cdot \Phi(z)| = +\infty$, which indicates that \,$(G_n \cdot \Phi)$ \,does not converge to \,$0$ \,uniformly on \,$\Omega_i$, hence neither on \,$\Omega$. Thus, the proof is finished in the case that (i) holds.

\vskip 3pt

Finally, we tackle the case in which (ii) happens; that is, $N(\Omega)=\N$, there exist $\alpha>0$, a sequence \,$(z_n)$ \,of pairwise different points in \,$\Omega$, a subsequence \,$(F_n)$ \,of \,$(f_n)$ \,with \,$|F_n(z_n)|\geq \alpha$ \,for all \,$n\in\N$, and pairwise different indexes
\,$i(1),i(2),i(3),\ldots$ \,in \,$N(\Omega)$ \,such that \,$z_n\in \Omega_{i(n)}$ \,for all \,$n\in\N$. Let us split \,$\N$ \,into infinitely many pairwise disjoint infinite sets $E_1,E_2,E_3,\ldots$ and, for each \,$k\in\N$, consider consider the set
\begin{equation*}
  S_k := \begin{cases}
           \bigcup_{n\in E_k}\Omega_{i(n)}           & \mbox{if } k\geq 2 \vspace{2mm} \\
           \Omega\setminus \bigcup_{k=2}^{\infty}S_k & \mbox{if } k=1.
         \end{cases}
\end{equation*}
Observe that the \,$S_k$'s are pairwise disjoint open subsets of \,$\Omega$ \,and \,$\Omega=\bigcup_{k=1}^{\infty}S_k$. Moreover, $z_n\in\Omega_{i(n)}\subset S_k$ \,for all \,$n\in E_k$, so for each \,$k\in\N$ \,there is an \,$m(k)\in\N$ \,satisfying \,$f_{m(k)}|_{S_k}\neq 0$. By Lemma \ref{Lemma para espaciabilidad de Sp-Suc y Suc-Su}(d), the vector subspace of \,$H(\Omega)^{\N}$ \,defined by
\begin{equation*}
  X := \left\{ \sum_{k=1}^{\infty} c_k \cdot \chi_{S_k} : c_k\in\C \hbox{ for all } k\in\N \right\}
\end{equation*}
is closed, infinite dimensional, and \,${\bf f} \in M({\bf f})\subset \mathcal{S}_{uc}$.

\vskip 3pt

To conclude the proof, fix \,${\bf g}=(g_n)\in M({\bf f}) \setminus \{0\}$. Our final task is to show that \,${\bf g} \not\in \mathcal{S}_u$. By the definition of \,$M({\bf f})$, there is \,$\Phi = \sum_{k=0}^{\infty} c_k \cdot \chi_{S_k} \in M({\bf f}) \setminus \{0\}$ \,such that \,$c_k\neq 0$ \,for some \,$k\in \N$ \,and \,$g_n = f_n \cdot \Phi$ \,for all \,$n\in\N$. Observe that
\,$E_k$ \,is an infinite set and if \,$n\in E_k$, then
\begin{equation*}
  \sup_{z\in\Omega} |F_n(z) \cdot \Phi(z)|
  \geq \sup_{z\in S_k} |F_n(z) \cdot \Phi(z)|
  \geq |F_n(z_n) \cdot \Phi(z_n)|
  \geq \alpha \cdot |c_k| > 0.
\end{equation*}
Thus, $\sup_{z\in\Omega} |F_n(z) \cdot \Phi(z)|\geq \alpha \cdot |c_k| > 0$ \,for infinite many values of \,$n$. Therefore, we deduce that $(F_n \cdot \Phi)\not\in  \mathcal{S}_u$ and, thus, ${\bf g} \not\in \mathcal{S}_u$, as required.
\end{proof}

\section{Final remarks}

\noindent {\bf 1.} If $\Omega$ is a domain, the space $X$ constructed in the proof of Theorem \ref{Thm-Sp-Suc spaceable} is simply \,$H(\Omega)$, because \,$H(\Omega) = \langle 1 \rangle \oplus \left\{f\in H(\Omega): f(a)=0\right\}$ (where \,$a$ \,is any fixed point in \,$\Omega$). The reason why we cannot also select as \,$X$ \,the whole space \,$H(\Omega)$ \,in the proof of the Theorem \ref{Thm-Suc-Su spaceable} is that, due to the Weierstrass interpolation theorem (see, e.g., \cite[Theorem 15.11]{Rudin}), we can find a nonzero \,$\Phi \in H(\Omega)$ \,with zeros at each point \,$z_n$, which counteracts the relative largeness of the members of the subsequence \,$(F_n)$ \,at such points.

\vskip 5pt

\noindent {\bf 2.} Another interesting mode of convergence is the weak one, that is associated to the weak topology \,$\tau_w$ \,in \,$H(\Omega)$. This topology is strictly weaker than the natural one in \,$H(\Omega)$, and a basis of $0$-neighbourhoods on \,$\tau_w$ \,consists of the sets \,$\left\{f\in H(\Omega): |\Lambda (f)| < \varepsilon\right\}$, where \,$\varepsilon>0$ \,and \,$\Lambda$ \,runs over the topological dual \,$H(\Omega)^*$ \,of \,$H(\Omega )$. For a description of \,$H(\Omega)^*$ see, e.g., \cite{Grothendieck} or \cite{KotheDualitat}. The topology \,$\tau_w$ \,is non-metrizable, and a net $(f_\alpha)\subset H(\Omega)$ \,$\tau_w$-converges to \,$f$ \,if and only if \,$\Lambda(f_\alpha) \to \Lambda (f)$ \,for all \,$\Lambda \in H(\Omega)^*$.

Nevertheless, for {\it sequences} \,$(f_n)$, compact convergence and weak convergence are {\it equivalent.} Indeed, if \,$f_n \to f$ \,weakly, the set \,$\left\{f_n : n\in \N\right\}$ \,is weakly bounded. But \,$H(\Omega)$ \,is a locally convex space, which entails that \,$\{f_n : \, n \in \N \}$ \,is also bounded for the natural topology of \,$H(\Omega)$ (see, e.g., \cite[Theorem 3.18]{RudinFunctionalAnalysis}). This means that \,$\left\{f_n : n \in \N\right\}$ \,is uniformly bounded on each compact subset of \,$\Omega$. Since the point evaluations are in \,$H(\Omega)^*$, we get that \,$f_n(z) \to f(z)$ \,for all \,$z \in \Omega$. At this point, Vitali's convergence theorem tells us that \,$f_n \to f$ \,compactly on \,$\Omega$ (see, e.g., \cite[p.~154]{Conway}). Consequently,
\begin{equation*}
  \mathcal{S}_{uc} = \mathcal{S}_w :=
  \left\{(f_n) \in H(\Omega)^{\N} : f_n \to 0 \hbox{\, weakly} \right\},
\end{equation*}
which makes it unnecessary a further study of comparison between \,$\mathcal{S}_p$, $\mathcal{S}_{uc}$, $\mathcal{S}_{u}$, and \,$\mathcal{S}_{w}$.

\vskip 6pt
%
%

\noindent {\bf Authors contribution.} All four authors have contributed equally to this work.

\end{document}